\documentclass[12pt]{amsart}
\usepackage{amsmath, geometry, amssymb} 
\numberwithin{equation}{section}
\newtheorem{theorem}{Theorem}[section]
\newtheorem{lemma}{Lemma}[section]
\newtheorem{corollary}{Corollary}[section]
\newtheorem{proposition}{Proposition}[section]

\newtheorem{notation}{Notation}[section]

\title{Shift-Inequivalent Decimations of the Sidelnikov-Lempel-Cohn-Eastman Sequences}
\author{\c{S}aban Alaca and Goldwyn Millar}
\date{} 

\begin{document}

\maketitle

\begin{abstract} 
We consider the problem of finding maximal sets of shift-inequivalent decimations of Sidelnikov-Lempel-Cohn-Eastman (SLCE) sequences (as well as the equivalent problem of determining the multiplier groups of the almost difference sets associated with these sequences). This is an open problem that was originally posed in \cite{C3} and that was mentioned more recently as being open in \cite{A4}.\\
We derive a numerical necessary condition for a residue to be a multiplier of an SLCE almost difference set. Using our necessary condition, we show that if $p$ is an odd prime and $S$ is an SLCE almost difference set over $\mathbb{F}_p,$ then the multiplier group of $S$ is trivial. Consequently, for each odd prime $p,$ we obtain a family of $\phi(p-1)$ shift-inequivalent balanced periodic sequences (where $\phi$ is the Euler-Totient function) each having period $p-1$ and nearly perfect autocorrelation. \\


\end{abstract}

\section{Introduction}

Let $\mathbf{a} = a_0 a_1 a_2 \ldots$ be a sequence of elements from the ring $\mathbb{Z}/M\mathbb{Z},$ where $M$ is some positive integer. 
Then $\mathbf{a}$ is \emph{periodic}
if there is an integer $v>0$ such that  $a_i = a_{v+i}$ for all integers $i\geq 0$.
If $v$ is the smallest such integer, then we say that $\mathbf{a}$ has period $v$; indeed, for the rest of this section, assume $\mathbf{a}$ is in fact periodic of period $v$. We shall discuss families of periodic sequences with certain special properties that have a variety of applications (for instance, in stream-cipher cryptography and code division multiple access (CDMA) communications systems). 

Let $\mathbf{b} = b_0 b_1 b_2 \ldots$ be another sequence of elements from the ring $\mathbb{Z}/M\mathbb{Z}.$ Assume that $\mathbf{b}$ is also periodic of period $v.$ The (periodic) correlation $\mathcal{C}_{\mathbf{a},\mathbf{b}}$ of $\mathbf{a}$ and $\mathbf{b}$ is defined as follows: for each nonnegative integer $\tau,$ \[\mathcal{C}_{\mathbf{a},\mathbf{b}}(\tau) := \sum_{t = 0}^{v-1} \text{exp}\left({\frac{2\pi i(a_t - b_{t + \tau})}{M}}\right),\]
where the terms of $\mathbf{a}$ and $\mathbf{b}$ appearing in the exponents are interpreted as integers. 

The function $\mathcal{C}_{\mathbf{a},\mathbf{a}}$ is called the \emph{autocorrelation} of $\mathbf{a},$ and the values $\mathcal{C}_{\mathbf{a},\mathbf{a}}(\tau)$ for $1\leq \tau \leq v-1$ are called the \emph{out-of-phase} autocorrelation values of $\mathbf{a}.$ We say that $\mathbf{a}$ has \emph{low out-of-phase autocorrelation} if its out-of-phase autocorrelation values are small compared to $v.$ Low out-of-phase autocorrelation is one of the criteria for a periodic sequence to be suitable for use as a key sequence in a stream-cipher cryptosystem (see the discussion of Golomb's postulate R3 in \cite[Section 5.1]{G1} and \cite{W1}). Indeed, for a sequence to be useful for this purpose, its out-of-phase autocorrelation values should be close to zero.

If there exists an integer $\ell$ such that for each positive integer $i,$ $a_i = b_{i+\ell},$ then we say that $\mathbf{a}$ and $\mathbf{b}$ are \emph{shift-equivalent} and that $\mathbf{a}$ and $\mathbf{b}$ are \emph{shifts} of one another; in this case, $\mathcal{C}_{\mathbf{a},\mathbf{b}}$ can be obtained from $\mathcal{C}_{\mathbf{a},\mathbf{a}}$ by simple formulae. If $\tau \geq \ell,$ then  $C_{\mathbf{a},\mathbf{b}}(\tau) = C_{\mathbf{a},\mathbf{a}}(\tau - \ell)$; if $\tau < \ell,$ then  $C_{\mathbf{a},\mathbf{b}}(\tau) = C_{\mathbf{a},\mathbf{a}}(\tau + v - \ell).$ If $\mathbf{a}$ and $\mathbf{b}$ are shift inequivalent, then we say that $\mathcal{C}_{\mathbf{a},\mathbf{b}}$ is the \emph{cross-correlation} of $\mathbf{a}$ and $\mathbf{b}.$ Furthermore, we consider a family $\mathcal{F}$ of shift-inequivalent sequences to have \emph{low cross-correlation} if for any pair of sequences $\mathbf{c},\mathbf{d}$ in $\mathcal{F}$ and for any $\tau,$ $\mathcal{C}_{\mathbf{c},\mathbf{d}}(\tau)$ is small compared to $v.$ More precisely, following \cite{G1}, we stipulate that $\mathcal{F}$ has low cross-correlation if for any pair of sequences $\mathbf{c},\mathbf{d}$ from $\mathcal{F},$ $\mathcal{C}_{\mathbf{c},\mathbf{d}}$ outputs only values less than or equal to $\delta\sqrt{v} + \epsilon,$ for some small integers $\delta$ and $\epsilon.$
 
For CDMA applications, one would like to have families of shift-inequivalent sequences with low cross-correlation such that each member of each family has low out-of-phase autocorrelation (see, for instance, \cite{G1} or \cite{G3}). Furthermore, one would like such families to be as large as possible (i.e. to include as many sequences as possible). One may also desire that sequences in these families have certain additional properties (such as cryptographic strength). 

One indicator of cryptographic strength (besides low out-of-phase autocorrelation) is the \emph{balance property}. We say that the sequence $\mathbf{a}$ is \emph{balanced} if in a given period of $\mathbf{a}$ (i.e. in a given list of $v$ consecutive elements of $\mathbf{a}$) each element of $\mathbb{Z}/M\mathbb{Z}$ appears either $\lfloor v/M \rfloor$ or $\lceil v/M \rceil$ times.

There are several known families of shift-inequivalent sequences with good periodic correlation properties. We begin by briefly discussing the families most relevant to our work in this paper. 

Let $p$ be a prime, let $d$ be a positive integer, and let $q = p^d.$ Let $\alpha$ be a primitive element of $\mathbb{F}_q,$ let $A \in \mathbb{F}_q,$ and let $\text{Tr}$ denote the field trace from $\mathbb{F}_q$ to $\mathbb{F}_p$ given by the rule that for $\beta \in \mathbb{F}_q.$ $\text{Tr}(\beta) = \beta + \beta^p + \cdot \cdot \cdot + \beta^{p^{d-1}}.$ Let $\mathbf{m} = m_0 m_1 m_2 \ldots$ be the sequence defined by the rule that $m_n = \text{Tr}(A\alpha^{-n}).$ Then $\mathbf{m}$ is called an \emph{m-sequence} of degree $d$ over $\mathbb{F}_p.$ 

The basic properties of m-sequences were discovered in the 1950s by Golomb \cite{G4} and Zierler \cite{Z1}. In the case that $p = 2,$ an m-sequence is (essentially) a combinatorial object called a \emph{Singer difference set}. These objects were originally discovered and studied by Singer in the 1930s \cite{S1}. 

It is known that the m-sequence $\mathbf{m}$ is a balanced sequence with period $q-1.$ The m-sequences also have near ideal autocorrelation: if $\tau \neq 0,$ then $\mathcal{C}_{\mathbf{m},\mathbf{m}}(\tau) = -1.$ Furthermore, these sequences have another desirable cryptographic property called the \emph{run} property (see \cite[Chapter 5]{G1} or \cite[Chapter 10]{G3} for proofs of all of these claims). However, the m-sequences do have \emph{low linear-complexity}, which is a type of cryptographic weakness (see \cite[Section 5.1]{G1}).

It is possible to use m-sequences to build families of shift-inequivalent, balanced, periodic sequences with good correlation properties. Let $t \geq 1$ be an integer. Then the \emph{t-fold decimation} of $\mathbf{a}$ is the sequence whose $i$th entry is $a_{ti}.$ Following \cite{G3}, we denote the t-fold decimation of $\mathbf{a}$ by $\mathbf{a}[t].$ For a proof of the following result, see \cite[Proposition 10.2.1]{G3}.
\begin{lemma} \label{m-lem} Let $p$ be a prime, let $d$ be a positive integer, and let $q = p^d.$ Let $\mathbf{m}$ be an m-sequence of degree $d$ over $\mathbf{F}_p.$\\
$1)$ Every m-sequence of degree $d$ over $\mathbf{F}_p$ is a shift of a decimation of $\mathbf{m}.$\\
$2)$ The decimation $\mathbf{m}[t]$ is again an m-sequence if and only if $t$ is relatively prime to $q-1.$\\
$3)$ The decimation $\mathbf{m}[t]$ is a shift of $\mathbf{m}$ if and only if $t$ is a power of $p.$\\
$4)$ There are $\phi(q-1)/d$ shift inequivalent m-sequences of degree $d$ over $\mathbb{F}_p.$
\end{lemma}

Thus, the set of decimations $\lbrace \mathbf{m}[t] \rbrace,$ where $t$ ranges over a set of integers congruent to representatives of the distinct cosets of $\langle p \rangle$ in $\mathbb{Z}/(q-1)\mathbb{Z}^*$, is a family of $\phi(q-1)/d$ shift inequivalent m-sequences. It is still an open problem to determine the precise cross-correlation values of the pairs of sequences in this family. However, cross-correlation values are known in certain special cases. If $t = 1 + p^i$ for some $i,$ then we say that $\mathbf{m}[t]$ is a \emph{quadratic decimation}, and the precise values taken on by the function $\mathcal{C}_{\mathbf{m},\mathbf{m}[t]}$ are known (see \cite{G2}, \cite{H1}, \cite{M1}, and \cite{P1}; alternatively, see the discussion in \cite{G3}). Also, if $t = -1,$ then in some cases, the values taken on by $\mathcal{C}_{\mathbf{m},\mathbf{m}[t]}$ are known (see \cite{L2}; the results from \cite{L2} are also briefly discussed in \cite{G3}).

We define the \emph{termwise sum} of the sequences $\mathbf{a}$ and $\mathbf{b}$ to be the sequence whose $i$th term is $a_i + b_i.$ Several authors have considered families of termwise sums of m-sequences having the same period. One such family is the family of Gold sequences, which is a family comprised of sequences constructed by taking termwise sums of m-sequences with shifts of their quadratic decimations \cite{G8}. The family of Gold sequences of a given period is larger than the family of shift-inequivalent decimations of an m-sequence having the same period, and it also has nice cross-correlation properties. However, the Gold sequences have worse autocorrelation properties than the m-sequences (and they are not always balanced). The correlation properties of the Gold sequences are summarized in \cite[Table 11.1]{G3}. Interestingly, the Gold sequences are currently used in the civilian C/A code for the US GPS system (see \cite[Section 11.2, Exercise 2]{G3}). 

Finally, m-sequences have been used to construct another class of sequences called Gordon-Mills-Welch (GMW) sequences. These sequences were first discovered (in the binary case, i.e. the case in which $p = 2$) by Gordon, Mills, and Welch \cite{G6}. The construction of the GMW sequences  relies on some rather deep results concerning the ``array structure'' of m-sequences. GMW sequences are discussed in \cite{G1} and \cite{G3}. For a nice discussion of the binary version of these sequences, see \cite{B1}.

In this paper, we consider another class of sequences, which are similar to the m-sequences in that their definition relies on both the multiplicative and additive structures of finite fields. Let $p$ be an odd prime, let $d$ be a positive integer, and let $q = p^d.$ Let $\alpha$ be a primitive element of $\mathbb{F}_q,$ and let $M|q-1.$ Following \cite{G5}, for $0 \leq k \leq M-1,$ we set $D_k = \lbrace \alpha^{Mi + k}-1| 0 \leq i < (q-1)/M \rbrace.$ An M-ary Sidelnikov sequence $\mathbf{s} = s_0 s_1 s_2 \ldots$ is a sequence of period $q-1$ whose first $q-1$ elements are defined as follows: for $0 \leq j < q-1,$ \[s_j = \begin{cases} 0 & \text{ if } \alpha^j = -1\\ k & \text{ if } \alpha^j \in D_k \end{cases}.\] 

This class of sequences was originally discovered by Sidelnikov in 1969 \cite{S2} and, in the binary case (i.e. the case in which $M = 2$) rediscovered independently by Lempel, Cohn, and Eastman in 1977 \cite{C3}. Consequently, we follow the authors of \cite{K1} in referring to the M-ary Sidelnikov sequences as Sidelnikov-Lempel-Cohn-Eastman (SLCE) sequences in the case that $M = 2.$

The Sidelnikov sequences have low out-of-phase autocorrelation. Indeed, in the case that $M = 2,$ if $\frac{1}{2}(q-1)$ is odd, then every out-of-phase autocorrelation value of $\mathbf{s}$ is either $\pm 2,$ and if $\frac{1}{2}(q-1)$ is even, then every out-of-phase autocorrelation value of $\mathbf{s}$ is either $0$ or $-4$ (see \cite{C3}). It is also clear from the definition of these sequences that they have the balance property.

If $c \in \mathbb{Z}/M\mathbb{Z},$ then we stipulate that $c\mathbf{a}$ is the sequence whose $i$th entry is $ca_i$ and we say that $c\mathbf{a}$ is a \emph{constant multiple} of $\mathbf{a}.$ The authors of \cite{K4} use the Weil bound (a version of which is given in Theorem \ref{wb} of the present paper) to prove an upper bound on the cross-correlation of a two distinct constant multiples of a Sidelnikov sequence.
\begin{theorem} \label{good bound} \cite{K4} Let $q$ be a power of an odd prime, and let $M|q-1.$ Let $\mathbf{s}$ be an M-ary Sidelnikov sequence over $\mathbb{F}_q^*.$ Let $c_1,c_2 \in \mathbb{Z}/M\mathbb{Z},$ $c_1,c_2 \neq 0.$ Let $\mathbf{a} := c_1\mathbf{s},$ and let $\mathbf{b} := c_2 \mathbf{s}.$ The for each $\tau = 0,...,q-2,$ \[|\mathcal{C}_{\mathbf{a},\mathbf{b}}(\tau)| \leq \sqrt{q} + 3.\]
\end{theorem}
As the authors of \cite{K4} note, it follows from Theorem \ref{good bound} that the set of all nonzero constant multiples of an $M$-ary Sidelnikov sequence forms a set of $M-1$ sequences with nearly ideal cross-correlation.

Using the Weil bound, the authors of \cite{K3} derive an upper bound on the cross-correlation of two shift-inequivalent decimations of constant multiples of a Sidelnikov sequence.
\begin{theorem} \label{bound from weil} \cite{K3} Let $q$ be a power of an odd prime, let $d,d^{\prime} \in \mathbb{Z},$ and let $M|q-1.$ Let $\mathbf{s}$ be an M-ary Sidelnikov sequence over $\mathbb{F}_q^*.$ Assume that $(d,q-1) = (d^{\prime},q-1) = 1$ and that $p$ divides neither $d$ nor $d^{\prime}.$ Let $c_1,c_2 \in \mathbb{Z}/M\mathbb{Z},$ $c_1,c_2 \neq 0.$ Let $\mathbf{a} := c_1\mathbf{s}[d],$ and let $\mathbf{b} := c_2 \mathbf{s}[d^{\prime}].$ Assume that $\mathbf{a}$ and $\mathbf{b}$ are shift-inequivalent. Then for each $\tau = 0,...,q-2,$ \[|\mathcal{C}_{\mathbf{a},\mathbf{b}}(\tau)| \leq (d+d^{\prime}-1)\sqrt{q} + 3.\]
\end{theorem}
The authors of \cite{K3} explicitly computed the cross-correlations of two shift-inequivalent decimations of constant multiples of a Sidelnikov sequence in several particular cases. In each case they considered, they found that the actual cross-correlation values were well below the upper bounds implied by Theorem \ref{bound from weil}.
 
Several authors have used Sidelnikov sequences to construct families of sequences with low cross-correlation in a manner similar to the way that m-sequences are used to construct the Gold sequences. 
In \cite{C1}, the authors consider a family consisting of term-wise sums of constant multiples of a Sidelnikov sequence with constant multiplies of one of its shifts. Let $\mathbf{s}$ be an $m$-ary Sidelnikov sequence over $\mathbb{F}_q^*,$ where $q = p^d.$ For $1 \leq c_1, c_2 \leq M-1$ and $0 \leq r \leq q-1,$ let $u_{c_1,c_2;r}$ be the sequence whose $i$th entry $u_{c_1,c_2;r}(i)$ is defined by $u_{{c_1,c_2};r}(i) := c_1s_{i} + c_2 s_{i+r}.$ Let $T:= \lceil \frac{q-1}{2} \rceil.$  Let 
\begin{eqnarray*} \mathcal{L} := \lbrace u_{{c_1,0};0} | 1 \leq c_1 \leq M-1 \rbrace \\
\cup \lbrace u_{c_1,c_2;i} | 1 \leq c_1,c_2 \leq M-1, 1 \leq i \leq T-1 \rbrace \\
\cup \lbrace u_{c_1,c_2;T} | 1 \leq c_1 < c_2 \leq M-1 \rbrace
\end{eqnarray*}
The authors of \cite{C1} use the Weil bound to obtain an upper bound on the cross-correlation values of the sequences in $\mathcal{L}.$
\begin{theorem} \cite{C1} \label{prelim upbd} The family $\mathcal{L}$ consists of $(M-1)^2(T-1) + M(M-1)/2$ shift-inequivalent sequences. Furthermore, the magnitudes of the cross-correlation values of any two distinct sequences in $\mathcal{L}$ are less than or equal to $3\sqrt{q}+5.$
\end{theorem} 

The authors of \cite{C2} enlarge the family from \cite{C1} by adding in termwise sums of constant multiples of a Sidelnikov sequence with constant multiples of its decimation by $-1.$ For $1 \leq c_1, c_2 \leq M-1$ and $0 \leq r < q-1,$ let $v_{c_1,c_2;r}$ be the sequence whose $i$th entry $v_{c_1,c_2;r}(i)$ is defined by $v_{{c_1,c_2};r}(i) := c_1s_{i} + c_2 s_{-i+r}.$ Let 
\begin{eqnarray*} \mathcal{K} := \lbrace v_{{0,c_1};0} | 1 \leq c_1 \leq M-1 \rbrace \\
\cup \lbrace v_{c_1,c_2;i} | 1 \leq c_1,c_2 \leq M-1, 1 \leq i \leq T-1 \rbrace \\
\cup \lbrace v_{c_1,c_2;T} | 1 \leq c_1,c_2 \leq M-1, c_1 \neq c_2 \rbrace
\end{eqnarray*}
Let $\mathcal{M} := \mathcal{K} \cup \mathcal{L}.$ The authors of \cite{C2} use the Weil bound obtain an upper bound on the cross-correlation values of the sequences in $\mathcal{M}.$
\begin{theorem} \cite{C2} \label{upbd} The family $\mathcal{M}$ consists of $2(M-1)^2(T-1) + 2(M-1) + 3(M-1)(M-2)/2$ shift-inequivalent sequences. Furthermore, the magnitudes of the cross-correlation values of any two distinct sequences from $\mathcal{M}$ are less than or equal to $4\sqrt{q} + 5.$
\end{theorem}
Thus, the family $\mathcal{M}$ is nearly twice as large as the family $\mathcal{L}.$ However, the cross-correlation of $\mathcal{M}$ is slightly worse than the cross-correlation of $\mathcal{L}.$   

We conclude our preliminary discussion of the Sidelnikov sequences by noting that the authors of \cite{G5} show that certain Sidelnikov sequences have a nice ``array structure'' somewhat analogous to the ``array structure'' of the m-sequences which is used in the construction of the GMW sequences. Furthermore, they make use of this ``array structure'' to generate a family of sequences with good correlation properties, and they show that this family can be combined with the family from \cite{C1} to form an even larger family with good correlation properties. We note also that the authors of \cite{K2} have extended the results from \cite{G5}. 

In this paper, we consider the problem of determining the shift-inequivalent decimations of the SLCE sequences (i.e. the binary Sidelnikov sequences). Thus, we attempt to obtain an analogue of Lemma \ref{m-lem} (in the binary case). 

To that end, we prove a result that gives an easily checkable sufficient condition to determine whether two decimations of an SLCE sequence are shift-inequivalent. Using our result, we are able to show that an SLCE sequence over $\mathbb{F}_p^*$ is always shift-inequivalent to each of its decimations. Consequently, we are able to produce families of shift-inequivalent sequences with good autocorrelation properties. 

Instead of studying SLCE sequences directly, it is convenient for us to instead consider certain  combinatorial objects that are closely related to periodic binary sequences with good autocorrelation. Let $v \in \mathbb{Z}.$ If $A \subset \mathbb{Z}/v\mathbb{Z},$ we say that the \emph{characteristic sequence} of $A$ is the sequence $\mathbf{a}$ over $\mathbb{F}_2$ of period $v$ such that for $i = 0,...,v-1,$ $a_i = 1$ if $i \in A$ and $a_i = 0$ otherwise. We also say that $\mathbf{a}$ is the sequence \emph{associated with} $A$ and that $A$ is the set of residues \emph{associated with} $\mathbf{a}.$

Let $v,$ $k.$ and $\lambda$ be positive integers. A set $D$ of $k$ residues mod $v$ is called a $(v,k,\lambda)$ \emph{cyclic difference set} if for each nonzero residue $x$ mod $v,$ there exist exactly $\lambda$ ordered pairs of elements $y_1,y_2 \in D$ for which $y_1 - y_2 = x$ (see, for instance, \cite[Section VI]{B3} for an overview of the theory of difference sets). The sequence $\mathbf{d}$ associated with $D$ has two-valued autocorrelation
\[ \mathcal{C}_{\mathbf{d},\mathbf{d}}(\tau) = \begin{cases} v-4(k-\lambda) & \text{ if } \tau \neq 0\\ v & \text{ if } \tau = 0 \end{cases}\]
(see \cite[Section 7.2]{G1}). 
As we mentioned at the bottom of p.2, the m-sequences over $\mathbb{F}_2$ are associated with a class of cyclic difference sets called \emph{Singer difference sets}.

Let $r$ be a positive integer. A set $E$ of $k$ residues mod $v$ is called a $(v,k,\lambda,r)$ \emph{cyclic almost difference set} if there exists a set $R$ of $r$ nonzero residues mod $v,$ each of which can be written as a difference of elements of $E$ in exactly $\lambda$ ways and if every other nonzero residue mod $v$ can be written as a difference of elements of $E$ in exactly $\lambda + 1$ ways (see \cite{A2} or \cite{N1} for surveys of the theory of almost difference sets). The sequence $\mathbf{e}$ associated with $E$ has three-valued autocorrelation 
\[ \mathcal{C}_{\mathbf{e},\mathbf{e}}(\tau) = \begin{cases} v-4(k-\lambda) & \text{ if } \tau \in R\\ v-4(k - \lambda - 1) & \text{ if } \tau \in \mathbb{Z}/v\mathbb{Z}-\lbrace R \cup \lbrace 0 \rbrace \rbrace\\ v & \text{ if } \tau = 0 \end{cases}\]
(see \cite[Theorem 2]{N1}).
As we mentioned on p.4, the SLCE sequences have three-valued autocorrelation. Hence, these sequences are associated with a class of cyclic almost difference sets. Indeed, let $p$ be an odd prime, let $d$ be a positive integer, and let $q = p^d.$ Let $\mathbf{s}$ be an SLCE sequence defined over $\mathbb{F}_q^*.$ We will refer to the cyclic almost difference set $S$ corresponding to the SLCE sequence $\mathbf{s}$ as an SLCE cyclic almost difference set.

We now introduce another way of thinking about sequences and cyclic difference sets/almost difference sets. Let $G$ be a finite cyclic group of order $v$. The integral group ring $\mathbb{Z}[G]$ consists of 
all formal sums $\sum_{g \in G} a_g g$, where $a_g \in \mathbb{Z}$ and with addition and multiplication defined as follows:
\begin{eqnarray*}
\sum_{g \in G} a_g g + \sum_{g \in G} b_g g = \sum_{g \in G} (a_g + b_g) g
\end{eqnarray*}
and
\begin{eqnarray*}
\Big( \sum_{g \in G} a_g g \Big) \Big(  \sum_{h \in G} b_h h \Big) 
= \sum_{f\in G} \Big( \sum_{gh =f} a_g  b_h \Big) f .
\end{eqnarray*}
For any subset $A \subseteq G,$ we identify $A$ with the group ring sum of all the elements in $A$; indeed, we refer to this sum as $A.$ 
\begin{notation} Throughout this paper, if $r$ is some element of $\mathbb{Z}/n\mathbb{Z}$ (where $n$ is a positive integer) then we will write $r^{\prime}$ to denote the least positive integer in the congruence class $r.$
\end{notation}
Let $t \in (\mathbb{Z}/v\mathbb{Z})^*$, let $A = \sum a_g g$ be a group ring element, and let $A^{(t)} := \sum a_g g^{t^{\prime}}.$ We say that $t$ is a \emph{multiplier} of $A$ if there exists $\alpha \in G$ such that $A^{(t)} = \alpha A,$ and we say that $t$ is a \emph{strong multiplier} of $A$ if $A^{(t)} = A.$ Let $A^c$ denote the complement of $A$ in $G.$ Note that $t$ is a multiplier of $A$ if and only if $t$ is a multiplier of $A^c.$ The elements of $(\mathbb{Z}/v\mathbb{Z})^*$ which are multipliers of $A$ form a group; we refer to this group as the \emph{multiplier group} $H$ of $A.$ Similarly, the elements of $(\mathbb{Z}/v\mathbb{Z})^*$ which are strong multipliers of $A$ also form a group; we refer to this group as the \emph{strong multiplier group} $H_0$ of $A.$ 

Let $E$ be a cyclic almost difference set, and let $\mathbf{e}$ be the sequence associated with $E.$ Then $E^{(t)}$ is the group ring element corresponding to the decimation $\mathbf{e}[t^\prime]$ of $\mathbf{e}.$ Furthermore, $t$ is a multiplier of $E$ if and only if $\mathbf{e}[t^\prime]$ is a shift of $\mathbf{e};$ likewise, $t$ is a strong multiplier of $E$ if and only if $\mathbf{e}[t^\prime] = \mathbf{e}.$ 

Let $I$ be a complete set of distinct coset representatives of $H.$ The set $\lbrace \mathbf{e}(t^{\prime}): t \in I\rbrace$ is a maximal set of shift-inequivalent decimations of $\mathbf{e}.$ Thus, the problem of determining the shift-inequivalent decimations of $\mathbf{e}$ is equivalent to the problem of determining the multiplier group $H.$ 

It is clear from Lemma \ref{m-lem} that when $p=2,$ $d \in \mathbb{N},$ and $q = p^d,$ the multiplier group of a Singer difference set corresponding to an m-sequence defined over $\mathbb{F}_q$ is $\langle 2 \rangle.$ Indeed, for any prime $p$ and any power $q=p^d,$ it is possible to define Singer difference sets over $\mathbb{F}_q,$ although they don't correspond to m-sequences when $p \neq 2$ (see \cite{S2}). It is known that the multiplier group of a Singer difference set defined over $\mathbb{F}_q$ is $\langle p \rangle$ (see \cite{G6}).

The problem of determining the multiplier groups of the SLCE cyclic almost difference sets was considered in \cite{C3}. The authors of \cite{C3} were able to show that $\langle p \rangle$ is a subgroup of the group $H_0$ of strong multipliers of $S.$ Furthermore, they explicitly computed the multiplier groups of SLCE cyclic almost difference sets in a number of cases. They found that for most of the SLCE cyclic almost difference sets that they considered, $\langle p \rangle$ comprised the entire multiplier group. However, they did find a case in which $\langle p \rangle$ was actually a proper subgroup of the multiplier group: when $p=3$ and $d=2,$ the multiplier group $H$ of the SLCE almost difference set $S$ is $(\mathbb{Z}/(3^2-1)\mathbb{Z})^* = \lbrace 1,3,5,7 \rbrace \neq \lbrace 1,3 \rbrace = \langle p \rangle.$ It is mentioned in \cite{A2} that the problem of determining the multiplier groups of the SLCE cyclic almost difference sets is still open. In this paper, we make some progress towards solving this problem using group characters and facts about cyclotomic fields.

Our approach is essentially a version the character method for studying difference sets (see \cite{B3}) and in particular is somewhat akin to the methods employed in \cite{E1}, \cite{X1}, and \cite{Y1}. However, some interesting complications arise from the fact that the objects we are studying are almost difference sets rather than difference sets.

Whereas the character values of the difference sets considered in \cite{E1} can be expressed as multiples of Gauss or Jacobi sums, it turns out that the character values of the SLCE sequences are multiples of Jacobi sums plus an extra constant term (see Lemma \ref{char} in Section 3). Thus, whereas Stickleberger's theorem can be brought directly to bear on problems concerning the difference sets considered in \cite{E1}, it is a little more difficult to apply Stickleberger's theorem to problems concerning the SLCE almost difference sets. We are able to surmount this difficulty by making use of a technical lemma about the norms of certain elements in cyclotomic fields (Lemma \ref{norm}). This lemma enables us to apply Stickleberger's theorem to our problem and thus to prove the main theorem of this paper (Theorem \ref{nec}).

We had several motivations for considering the problem of determining the multiplier groups of the SLCE almost difference sets.

As we have mentioned, knowing multiplier groups of SLCE almost difference sets allows one to produce maximal families of shift-inequivalent decimations of SLCE sequences. The relation of these families to the other Sidelnikov families is similar to the relation between the families of decimations of m-sequences and the Gold sequences. The families of decimations of SLCE sequences have fewer sequences than the other Sidelnikov families, but they do have the advantage that each of their sequences has nearly perfect autocorrelation.

Alternatively, if one knows that a decimation $\mathbf{s}[t^{\prime}]$ of an SLCE sequence $\mathbf{s}$ is shift-inequivalent to $\mathbf{s},$ then one can also use $\mathbf{s}$ and $\mathbf{s}[t^{\prime}]$ to construct a family of shift-inequivalent sequences using a construction similar to the one given in \cite{C2}.

Since the SLCE sequences are balanced and have good autocorrelation properties, they are candidates for use as key sequences in stream cipher cryptosystems. However, the existence of a non-trivial multiplier for the group ring element associated with a sequence is a cryptographic weakness. Suppose the sequence $\mathbf{s}$ is used as a key sequence for a stream cipher. If $\mathbf{s}$ is shift-equivalent to its decimation $\mathbf{s}[t^{\prime}]$ and if an eavesdropper intercepts part of a cipher text message enciphered using $\mathbf{s}$, then she may be able to recover more of the cipher text by correlating the part she has received with its decimation by $t^{\prime}.$ 

Explicit determinations of multiplier groups of difference sets and related combinatorial objects are also of interest in theoretical combinatorics. Multiplier theorems (i.e. theorems guaranteeing the existence of multipliers for combinatorial structures of certain types and having certain parameters) are useful for proving nonexistence results about difference sets and related objects (see \cite[Sections VI.2 and VI.4]{B3}), and explicit determinations of multiplier groups of known structures can help point the way to more general theorems. Furthermore, computations of multiplier groups can be useful for proving different types of combinatorial results. For instance, one of the authors made use of the determination of the multiplier groups of the Singer difference sets given in \cite{G6} to prove a structural result about a related class of combinatorial objects called circulant weighing matrices \cite{M2}.

Finally, we believe the problem under consideration in this paper is an intrinsically interesting design-theoretic problem. The SLCE almost difference sets are important and useful mathematical objects. So, it is natural to try to determine their symmetries. 

\section{Preliminary Results}

The use of characters in the study of cyclic difference sets dates back to the work of Marshall hall in the 1940s (see \cite{H2}). As Beth et al. mention \cite[bottom of p.10]{B3}, this approach became standard after the paper of Turyn from 1965 \cite{T1}. The use of characters in the study of multipliers dates back to the work of Yamamoto in 1963 \cite{Y1}; this approach was also employed by Xiang in 1994 \cite{X1}. For an overview of the use of characters in the theory of cyclic difference sets, see \cite[Section VI.3 and Sections VI.13-VI.16]{B3}. In this paper, we make particular use of characters defined over multiplicative groups of finite fields. The study of characters over finite fields was initiated by Gauss and Jacobi, who considered character sums now known as \emph{Gauss sums} and \emph{Jacobi sums}; for a discussion of characters over finite fields, see \cite[Chapter 8]{I1}.

\begin{notation} 
Let $n$ be a positive integer. Henceforth, the symbol $\zeta_n$ will denote a primitive, complex $n$th root of unity. 
\end{notation} 

Let $G$ be a cyclic group of order $v.$ A \emph{group character} is a homomorphism $\chi:G \to \langle \zeta_{v} \rangle.$ Such a homomorphism can be extended by linearity to a map from $\mathbb{Z}[G]$ to $\mathbb{Z}[\zeta_{v}],$  the ring of integers of the cyclotomic field $\mathbb{Q}(\zeta_v)$ of order $v.$ 

The \emph{order} of the character $\chi$ is equal to the largest order of the complex roots of unity $\chi(g)$ (as $g$ ranges over $G$). It is known that for each $n|v,$ there exist exactly $\phi(n)$ characters defined on $G$ having order $n.$ In fact, the characters of a group $G$ of order $v$ themselves form a group of order $v,$ which we shall denote $\widehat{G},$ under the operation of pointwise multiplication (so, for $\psi, \chi \in \widehat{G}$ and for $g \in G,$ $\psi \cdot \chi (g) = \psi(g) \chi(g)$).

We will make use of the following result, commonly known as the \emph{inversion formula}, relating character values to group ring elements (see \cite[Lemma VI.3.5]{B3}).

\begin{lemma} \label{inv} Let $G$ be a cyclic group of order $v,$ and let \[A = \sum_{h \in G}a_h h\] be an element in $\mathbb{Z}[G].$ Then the coefficients of $A$ can be recovered from the character values of $A$ as follows. For $h \in G,$ 
\[a_h = \frac{1}{|G|} \sum_{\chi \in \widehat{G}} \chi(A)\chi(h^{-1}).\]
Consequently, if for $A,B \in \mathbb{Z}[G]$ $\chi(A) = \chi(B)$ for every $\chi \in \widehat{G},$ then $A = B.$ 
\end{lemma}

We adopt the following convention. 
For an integer $i \in \lbrace 0, \ldots ,p-1 \rbrace,$ we  refer to the corresponding element of $\mathbb{F}_p$ by italicizing $\mathit{i}.$

We note that if $\chi$ is a nontrivial character on $\mathbb{F}_q^*,$ it is common to extend $\chi$ to a map on $\mathbb{F}_q$ by setting $\chi(\mathit{0}) = 0.$ However, it is sometimes useful to define $\chi(\mathit{0})$ to be equal to something else. In this paper, we consider both characters $\chi$ on $\mathbb{F}_q$ for which $\chi(\mathit{0}) = 0$ and characters $\chi$ on $\mathbb{F}_q$ for which $\chi(\mathit{0}) = 1.$ 

It is possible to define a logarithm over $\mathbb{F}_q.$ Let $q$ be an odd prime power, and let $\alpha$ be a primitive element of $\mathbb{F}_q.$ For $x \in \mathbb{F}_q,$ we stipulate that 
\[\text{log}_{\alpha}(x) = \begin{cases} i & \text{if } x = \alpha^i, \hspace{0.05in} 0 \leq i \leq q-2 \\
0 & \text{if } x = \mathit{0}.\end{cases}\]
As Gong and Yu note (\cite{G5} and \cite{G7}) the $M$-ary Sidelnikov sequence $\mathbf{s}$ defined over $\mathbb{F}_q$ using $\alpha$ is completely determined by the congruences \begin{equation} \label{id1} s_i \equiv \text{log}_{\alpha}(\alpha^i+1) \pmod{M}, \hspace{01in} 0 \leq i \leq q-2.\end{equation}
Still following Gong and Yu, we define a multiplicative character $\psi_M$ of order $M$ on $\mathbb{F}_q$ by the rule that for $x \in \mathbb{F}_q,$ 
\begin{equation} \label{id2} \psi_M(x) = \text{exp}\left( \frac{2\pi \text{log}_{\alpha}(x) i}{M}\right).\end{equation}
Note that $\psi_M(\mathit{0}) = 1.$

Gong and Yu remark that (\ref{id1}) and (\ref{id2}) imply the following identity, which we will use later in this paper.
\begin{equation} \label{id3} \text{exp}\left(\frac{2\pi s_j i}{M}\right) = \psi_M(\alpha^j + 1), \hspace{0.1in} 0 \leq j \leq q-2. \end{equation}

We need to make use of a well-known character sum bound called the Weil bound. There are several different versions of this result in the literature. The version we use is essentially the version stated as Corollary 2.3 in \cite{W2}. However, we also make use of the refinement introduced by Gong and Yu \cite{G7} to adapt the bound to character sums for which the characters $\psi$ involved satisfy the condition $\psi(0) = 1.$ In \cite{G7}, Gong and Yu apply their refinement to a slightly different version of the Weil bound than the one stated as Corollary 2.3 in \cite{W2}. However, their logic applies equally well to the result from \cite{W2}.
\begin{theorem} (Weil Bound) \label{wb} Let $f_1(x),...,f_n(x)$ be monic, pairwise prime polynomials in $\mathbb{F}_q[x]$ whose largest square-free divisors have degrees $d_1,...,d_n,$ respectively. Let $\psi_1,...,\psi_n$ be non-trivial characters on $\mathbb{F}_q.$ Suppose that for each $i = 1,...,n,$ $\psi_i(\mathit{0}) = 1.$ Assume that for some $1 \leq i \leq n,$ the polynomial $f_i(x)$ is not of the form $g(x)^{\text{ord}(\psi_i)}$ in $\mathbb{F}_q[x],$ where $\text{ord}(\psi_i)$ is the smallest positive integer $d$ such that $\psi_i^d = 1.$ For each $i = 1,...,n,$ let $e_i$ be the number of distinct roots of $f_i(x)$ in $\mathbb{F}_q.$ Then for any $a_i \in \mathbb{F}_q,$ $i = 1,...,n,$ 
\[ |\sum_{x \in \mathbb{F}_q} \psi_1(a_1f_1(x)) \cdot \cdot \cdot \psi_n(a_nf_n(x))| \leq \left(\sum_{i = 1}^n d_i -1\right)\sqrt{q} + \sum_{i = 1}^n e_i.\]

\end{theorem}

In this paper, we consider the following type of character sum (which is, in fact, a type of Jacobi sum).

\begin{notation} Let $k|q-1,$ and let $\chi$ be a character of order $k$ on $\mathbb{F}_q$ for which $\chi(\mathit{0}) = 0.$ Then we set \[ K(\chi) = \chi(\mathit{4})\sum_{x \in \mathbb{F}_q} \chi(x)\chi(1-x).\]
\end{notation} 

We note the following result for later reference (see \cite{B2}, Theorem $2.1.8$).
\begin{lemma} \label{cong} Let $\chi$ be a character on $\mathbb{F}_q$ of order $k>1$ for which $\chi(\mathit{0}) = 0.$ Then 
\[K(\chi) \equiv -q \pmod{2(1 - \zeta_k)}.\]
\end{lemma}

For each $j \in (\mathbb{Z}/v\mathbb{Z})^*,$ let $\sigma_j$ denote the automorphism of $\mathbb{Q}(\zeta_v)$ that maps $\zeta_v$ to $\zeta_v^{j^{\prime}}.$ Recall that $\text{Gal}(\mathbb{Q}(\zeta_v)/\mathbb{Q}) = \lbrace \sigma_j|j \in (\mathbb{Z}/v\mathbb{Z})^*\rbrace$ (see, for example, \cite[Theorem 26, p. 596]{D1}). 

We will need the folllowing technical lemma concerning cyclotomic fields (see \cite{B2}, Theorem $2.1.9$).
\begin{lemma} \label{norm} Let $k>1.$ Then the norm of $1-\zeta_k$ in $\mathbb{Q}(\zeta_k)$ is given as follows: 
\[N(1-\zeta_k) = \begin{cases} \ell & \text{ if $k$ is a power of a prime $\ell$} \\ 1 & \text{ otherwise}. \end{cases}\]
\end{lemma}
 
Let $m$ and $v$ be positive integers greater than $1,$ and let $m|v.$ It follows from the Fundamental Theorem of Galois Theory (specifically, \cite[Theorem 14.14, (3)]{D1}) that the extension $\mathbb{Q}(\zeta_v) \supset \mathbb{Q}(\zeta_m)$ is Galois (say, with Galois group $H$). Furthermore, by \cite[Theorem 14.14, (4)]{D1}, $\text{Gal}(\mathbb{Q}(\zeta_v)/\mathbb{Q})/H \cong \text{Gal}(\mathbb{Q}(\zeta_m)/\mathbb{Q}).$ So, if $N_v$ denotes the norm in $\mathbb{Q}(\zeta_v)$ and $N_m$ denotes the norm in $\mathbb{Q}(\zeta_m),$ then 
\[N_v(1-\zeta_m) = \prod_{\sigma \in \text{Gal}(\mathbb{Q}(\zeta_v)/\mathbb{Q})} \sigma(1-\zeta_m)  = 
\prod_{\sigma \in \text{Gal}(\mathbb{Q}(\zeta_m)/\mathbb{Q})} \sigma(1-\zeta_m)^{\phi(v)/\phi(m)}\] \[= N_m(1 - \zeta_m)^{\phi(v)/\phi(m)}.\]
Thus, we obtain the following corollary of Lemma \ref{norm}.
\begin{corollary} \label{general norm} Let $v$ and $m$ be positive integers greater than $1,$ and let $m|v.$ Then the norm of $1-\zeta_m$ in $\mathbb{Q}(\zeta_v)$ is given as follows: 
\[N(1-\zeta_m) = \begin{cases} \ell^{\phi(v)/\phi(m)} & \text{ if $m$ is a power of a prime $\ell$} \\ 1 & \text{ otherwise}. \end{cases}\]
\end{corollary}

Recall that the rings of integers of cyclotomic fields (and, indeed, the rings of integers of any algebraic number fields) are Dedekind domains (see, for instance, \cite[Theorem 8.1.1]{A1}). Hence, every ideal in $\mathbb{Z}[\zeta_v]$ factors uniquely as a product of prime ideals. 

\begin{notation} For the rest of this section, let $k$ be a positive integer that is not divisible by $p$ such that $d$ is the order of $p \pmod{k}.$ Let $q = p^d.$ Furthermore, let $P$ be a prime ideal lying over $(p)$ in $\mathbb{Z}[\zeta_k].$
\end{notation}

\begin{notation} Let $T$ be a set of distinct coset representatives of $\langle p \rangle$ in $(\mathbb{Z}/k\mathbb{Z})^*.$ 
\end{notation}

For a proof of the next result, see \cite[Theorem 13.2.2]{I1}.

\begin{theorem} \label{facp}  
\[(p) = \prod_{j \in T}\sigma_{j}(P).\]
Additionally, $\mathbb{Z}[\zeta_k]/P \cong \mathbb{F}_q,$ and the multiplicative subgroup $\lbrace \zeta_k^i + P| 0 \leq i <  k \rbrace$ is comprised of all $k$ of the $k$th roots of unity in $\mathbb{Z}[\zeta_k]/P.$
\end{theorem}
   
Let $\chi_P:(\mathbb{Z}[\zeta_k]/P)^* \to \mathbb{C}$ be the function defined by the rule that for $\alpha + P \in \mathbb{Z}[\zeta_k]/P,$ $\chi(\alpha + P) = \zeta_k^i,$ where $\zeta_k^i$ is the unique power of $\zeta_k$ congruent to $\alpha^{(q-1)/k} \pmod{P}$ (see \cite[Proposition 14.2.1 on p. 204]{I1}; see also \cite[(11.2.1) on p. 344]{B2}). Set $\chi_P(0 + P) = 0.$ Then $\chi_P$ is called a \emph{Techimuller character}. Any character of order $k$ on $\mathbb{F}_q$ can be viewed as a Techimuller character (so long as it maps $0 + P$ to $0$) by identifying a generator $\alpha + P$ of $(\mathbb{Z}[\zeta_k]/P)^*$ with a generator $\gamma$ of $\mathbb{F}_q^*$ such that $\chi_P(\alpha + P) = \chi(\gamma).$

The next theorem, which gives the prime ideal factorization of $(K(\chi))$ in $\mathbb{Z}[\zeta_k],$ is a consequence of Stickleberger's Theorem on Gauss sums (indeed, this result follows by setting $m = n = 1$ in Corollary $11.2.4$ from \cite{B2} and using Theorem $11.2.9$ from \cite{B2} to rewrite the term appearing in the exponent). 

\begin{theorem} \label{Stickleberger} \[(K(\chi_P)) = \prod_{j \in T} \sigma_{j^{-1}} (P)^{d - \sum_{i = 0}^{d-1} \lbrace \lfloor \frac{2j^{\prime}p^i}{k} \rfloor - 2\lfloor \frac{j^{\prime}p^i}{k} \rfloor \rbrace.}\]
\end{theorem}

As Berndt et al. note \cite[p. 349, comment in the proof of Theorem 11.2.9]{B2}, the term $\lbrace \lfloor \frac{2j^{\prime}p^i}{k} \rfloor - 2\lfloor \frac{j^{\prime}p^i}{k} \rfloor \rbrace$ appearing in Theorem \ref{Stickleberger} equals $1$ or $0$ according to whether the remainder upon dividing $j^{\prime}p^i$ by $k$ is greater than $k/2$ or not.

\section{Character Values}
We begin this section by fixing some notation that we will use for the rest of the paper.
\begin{notation} Let $Y:= \lbrace y \in \mathbb{F}_{q}^* \mid  y = x(\mathit{1}-x)   \text{ for some }  x \in \mathbb{F}_{q}^*  \rbrace ,$ and let $Z: = Y^c.$
\end{notation}

Our work in this paper relies on the following lemma of Lempel, Cohn, and Eastman \cite{C3}. \begin{lemma} \label{LCE lem}  
Let $q$ be an odd prime power, and let $S$ be an SLCE almost difference set over $\mathbb{F}_q^*.$ $Z$ is a shift of $S$: in fact, $Z = \mathit{-4}^{-1}S,$ so that $S = \mathit{-4}Z$ and $S^c = \mathit{-4}Y.$
\end{lemma} 
A version of the following lemma appeared in our recent paper \cite{A3}. For convenience (and since we are stating the result in slightly different language) we provide a proof.
\begin{lemma} \label{char} Let $q$ be an odd prime power, and let $S$ be an SLCE almost difference set over $\mathbb{F}_q^*.$ Let $\chi$ be a character on $\mathbb{F}_q.$ Then 
\[\chi(S^c) = \frac{1}{2}\chi(\mathit{-1})(K(\chi) + 1).\]
\end{lemma}
\emph{Proof.} The reasoning in the next two sentences is taken from \cite[Theorem 2.14]{B2},   
where it serves a different purpose. Let $\gamma \in \mathbb{F}_q^*$ be fixed. 
An element $x \in \mathbb{F}_q^*$  satisfies the equation  $x(\mathit{1}-x) = \gamma$ if and only if it satisfies the equation $(\mathit{2}x-\mathit{1})^2 = \mathit{1}-\mathit{4}\gamma.$ 
Hence, the number of solutions of the equation $x(\mathit{1}-x) = \gamma$ in $F_q^*$  is $1 + \rho(\mathit{1}-\mathit{4}\gamma)$, where  $\rho$ denotes the (unique) quadratic character on $\mathbb{F}_q$ (with $\rho(\mathit{0})$ set equal to $0$).
It follows that every element of $\mathbb{F}_q^*$ is represented either twice or zero times in the form $x(\mathit{1}-x),$ save for $\mathit{4}^{-1},$ which is represented once. 

Making use of Lemma \ref{LCE lem}, we see that   
\begin{eqnarray*}
&\chi(\mathit{-1})K(\chi) &= \chi(\mathit{-4})\sum_{x \in \mathbb{F}_q^*} \chi(x)\chi(\mathit{1}-x) \\
&&= \chi(\mathit{-4})\sum_{x \in \mathbb{F}_q^*} \chi(x(\mathit{1}-x)) = \chi(\mathit{-4})\chi \Big( \sum_{x\in \mathbb{F}_q^*} x(\mathit{1}-x) \Big) \\
&& = \chi(\mathit{-4})\chi( 2Y - \mathit{4}^{-1}) = \chi(2S^c - (\mathit{-1})) = 2\chi(S^c) - \chi(\mathit{-1}).
\end{eqnarray*}
So, we deduce that 
\[\chi(S^c) = \frac{1}{2} \chi(\mathit{-1})(K(\chi) + 1). \qed \]
Note that it follows from Lemma \ref{cong} that the expression we have given in Lemma \ref{char} for $\chi(S^c)$ is indeed an algebraic integer.

\section{Multiplier Theorems}

We begin by showing how the results of Section $3$ can be used to recover a theorem of Lempel, Cohn, and Eastman \cite{C3}.
\begin{theorem} \cite{C3} Let $q$ be a power of an odd prime $p$, and let $S$ be an SLCE almost difference set over $\mathbb{F}_q^*.$ Then \label{LCE} $\langle p \rangle$ is a subgroup of the strong multiplier group of $S.$
\end{theorem}
\emph{Proof (new).} Let $\chi$ be a character on $\mathbb{F}_q^*,$ and let $i$ be a positive integer. Then, by Lemma \ref{char} and the Child's Binomial Theorem,  
\[\chi((S^c)^{(p^i)}) = \sigma_{p^i}(\chi(S^c)) = \frac{1}{2} \sigma_{p^i}(\chi(\mathit{-1})(K(\chi) + 1))\]
\[ = \frac{1}{2} \chi(\mathit{-1})^{p^i} (\sum_{x \in \mathbb{F}_q^*} \chi(x^{p^i})\chi((\mathit{1}-x)^{p^i}) + 1) = \frac{1}{2} \chi(\mathit{-1}) (\sum_{x \in \mathbb{F}_q^*} \chi(x^{p^i})\chi(\mathit{1}-x^{p^i}) + 1)\]
\[ = \frac{1}{2}\chi(\mathit{-1})(K(\chi) + 1) = \chi(S^c).\]
It follows by Lemma \ref{inv} that $(S^c)^{(p^i)} = S^c.$ Consequently, $S^{(p^i)} = S.$ $\qed$ \\

As a result of Theorem \ref{LCE}, the problem of determining the multiplier group of $S$ reduces to determining which elements of $T$ are multipliers of $S.$ To that end, we now establish a necessary condition for an element $t \in T$ to be a multiplier of $S.$

\begin{theorem} \label{nec} Let $p$ be an odd prime, let $d$ be a positive integer, and let $q = p^d.$ Let $S$ be an SLCE almost difference set over $\mathbb{F}_q.$ Let $T$ be a set of distinct coset representatives of $\langle p \rangle$ in $(\mathbb{Z}/(q-1)\mathbb{Z})^*,$ and let $t\in T$ be a multiplier of $S.$ Then the sets $S_0 = \lbrace j \in T: d - \sum_{i = 0}^{d-1} \lbrace \lfloor \frac{2(j^{-1})^{\prime}p^i}{q-1} \rfloor - 2\lfloor \frac{(j^{-1})^{\prime}p^i}{q-1} \rfloor \rbrace >0 \rbrace$ and $tS_0$ are either identical or disjoint.
\end{theorem}
\emph{Proof.} Notice that $t$ is also a multiplier of $D = 2(\mathit{-1})S^c$ (where $2 \in \mathbb{Z}$ and $\mathit{-1} \in \mathbb{F}_{q-1}^*$). So, there exists $g \in \mathbb{F}_{q}^*$ such that $D^{(t)} = gD.$ Let $P$ be a prime ideal lying over $p$ in $\mathbb{Z}[\zeta_{q-1}].$ Recall that the Techimuller character $\chi_P$ can be treated as a character on $\mathbb{F}_{q}^*.$ By Lemma \ref{char}, we have that 
\[\chi_P(D^{(t)}) = \chi_P(g)\chi_P(D) = \zeta\chi_P(D) = \zeta(K(\chi_P) + 1),\]
where $\zeta$ is some (not necessarily primitive) $(q-1)$th root of unity. But, we also have that 
\[\chi_P(D^{(t)}) = \sigma_t(\chi_P(D)) = \sigma_t(K(\chi_P) + 1) = \sigma_t(K(\chi_P)) + 1.\]
Consequently, 
\begin{eqnarray} \label{zeta} \sigma_t(K(\chi_P)) - \zeta K(\chi_P) = \zeta - 1.\end{eqnarray}
Now, assume there is a prime ideal $Q$ lying over $p$ that contains both $K(\chi_P)$ and $\sigma_t(K(\chi_P)).$ Then, by (\ref{zeta}), $\zeta - 1 \in Q.$ Note that $\zeta$ is a primitive $m$th root of unity for some $m$ dividing $q-1.$ 

By Corollary \ref{general norm}, if $m$ is a product of more than one prime, then $N(1-\zeta) = 1,$ and it follows that $\zeta-1$ is a unit. However, since $\zeta - 1 \in Q,$ this implies that $Q = \mathbb{Z}[\zeta_{q-1}],$ which is a contradiction. On the other hand, if $m = \ell^s,$ for some prime $\ell$ and some positive integer $s,$ then by Corollary \ref{general norm}, $N(1 - \zeta) = \ell^{\phi(q-1)/\phi(m)},$ and so $\ell^{\phi(q-1)/\phi(m)} \in Q.$ But, since $\ell|m|(q-1),$ we have that $\text{gcd}(p,\ell^{\phi(q-1)/\phi(m)}) = 1.$ So, since $p \in Q$ also, the Euclidean Algorithm implies that $1 \in Q.$ Once again, we get the contradiction that $Q = \mathbb{Z}[\zeta_{q-1}].$

The remaining possibility is that $m = 1,$ i.e. that $\zeta = 1.$ In this case, (\ref{zeta}) implies that $\sigma_t(K(\chi_P)) = K(\chi_P).$ So, if \emph{any} prime ideal lying over $K(\chi_P)$ also lies over $\sigma_t(K(\chi_P)),$ then \emph{every} prime ideal lying over $K(\chi_P)$ also lies over $\sigma_t(K(\chi_P)).$ Consequently, by Theorem \ref{Stickleberger}, $S_0$ and $tS_0$ are either identical or disjoint. $\qed$   \\

As an application of Theorem \ref{nec}, we shall prove that in the case $d=1,$ the multiplier group of $S$ is trivial. But first, we identify a few special cases that cannot be dealt with using Theorem \ref{nec}. 

It is known that $-1$ is never a multiplier of a nontrivial cyclic difference set \cite{B1}. Interestingly, in the case that $p = 3$ and $d = 2,$ $-1$ actually is a multiplier of $S.$ Naturally, it is of interest to determine when $S$ has $-1$ as a multiplier. However, at least in the case $d=1,$ it is easy to see that our condition \emph{never} rules out $-1$ as a multiplier of $S.$ 

Let $d = 1.$ In this case, the set $S_0$ has a particularly simple description: namely, 
\[S_0 = \lbrace j\in (\mathbb{Z}/(p-1)\mathbb{Z})^*: (j^{-1})^{\prime} < (p-1)/2 \rbrace.\] Hence, it is clear that $j \in S_0$ if and only if $(j^{-1})^{\prime} < (p-1)/2$ if and only if $(-j^{-1})^{\prime} > (p-1)/2$ if and only if $-j \notin S_0.$ So, in this case, $S_0$ and $-S_0$ are disjoint and so our condition does not rule out $-1$ as a multiplier.

Fortunately, we do have another tool at our disposal to help determine whether $-1$ is a multiplier of $S.$ For, if $-1$ is a multiplier of $S,$ then for some $\tau = 0,...,q-2,$ $\mathcal{C}_{\mathbf{s},\mathbf{s}[-1]}(\tau) = q-1.$ However, by Theorem \ref{upbd}, $\mathcal{C}_{\mathbf{s},\mathbf{s}[-1]}(\tau) \leq 4\sqrt{q} + 5.$ So, we must have that $q - 1 \leq 4\sqrt{q} + 5,$ i.e. that $q-4\sqrt{q} \leq 6.$ But, for $x>4,$ $x - 4\sqrt{x}$ is an increasing function, and for $q = 3^3,$ $q - 4\sqrt{q} > 6.$ So, if $q \geq 27,$ then $-1$ is not a multiplier of $S.$ It can be checked directly that for $q < 27,$ $-1$ is a multiplier of $S$ if and only if $q $ equals $3$ or $9$. Thus, the following proposition is a direct consequence of the cross-correlation bound from \cite{C2}.

\begin{proposition} \label{minus one} Let $q$ be an odd prime power, and let $S$ be an SLCE almost difference set defined over $\mathbb{F}_q^*.$ Then $-1$ is a multiplier of $S$ if and only if $q$ equals $3$ or $9$.
\end{proposition}

There are two other cases we cannot handle using Theorem \ref{nec}. Let $d = 1,$ and let $p$ be a prime congruent to $1 \pmod{4}.$ Then Theorem \ref{nec} does not rule out $(p-1)/2 \pm 1$ as multipliers of $S.$

Note that $((p-1)/2 - 1)^{-1} = ((p-1)/2 - 1).$ Now, $j \in S_0$ if and only if $(j^{-1})^{\prime} < (p-1)/2$ if and only if $(p-1)/2 - (j^{-1})^{\prime} < (p-1)/2$ if and only if (since each element of $(\mathbb{Z}/(p-1)\mathbb{Z})^*$ is odd) $(((p-1)/2 - 1)j^{-1})^{\prime} < (p-1)/2$ if and only if $((p-1)/2 - 1)j \in S_0.$ So, $S_0 = ((p-1)/2 - 1)S_0.$ Similarly, one can argue that $S_0$ and $((p-1)/2 + 1)S_0$ are disjoint. It follows that we cannot use Theorem \ref{nec} to show that $(p-1)/2 \pm 1$ are not multipliers of $S$.

So, in order to show that $(q-1)/2 \pm 1$ are (almost) never multipliers of $S,$ we prove a (weak) bound on the cross-correlations of $\mathbf{s}$ with $\mathbf{s}[(q-1)/2 \pm 1]$ and make an argument similar to the one we made to show that $-1$ is (almost) never a multiplier of $S.$

\begin{lemma} \label{weakbd} Let $q$ be a prime power such that $q \equiv 1 \pmod{4},$ let $\alpha$ be a primitive element of $\mathbb{F}_q^*,$ and let $\mathbf{s}$ be the SLCE sequence defined over $\mathbb{F}_q^*$ using $\alpha.$ Then the cross-correlation values of $\mathbf{s}$ with $\mathbf{s}[(q-1)/2 \pm 1]$ are less than or equal to $\text{max}\lbrace 3\sqrt{q}+6, \frac{1}{2}(q + 3\sqrt{q} + 7)\rbrace.$
\end{lemma}
\emph{Proof.} In what follows, let $\rho$ be the quadratic character on $\mathbb{F}_q$, but let us stipulate that $\rho(\mathit{0}) = 1.$ Let $\tau = 0,...,q-2.$ Then, by (\ref{id3}), 
\[ |\mathcal{C}_{\mathbf{s}, \mathbf{s}[(q-1)/2 - 1]}(\tau)| 
= |\sum_{j = 0}^{q-2} \text{exp}\left(\frac{2\pi i (s_j - s_{((q-1)/2-1)j+\tau})}{2}\right)| \]
\[= |\sum_{j = 0}^{q-2} \rho(\alpha^j + 1) \rho(\alpha^{((q-1)/2 - 1)j + \tau} + 1)|\] \[ = |\sum_{\substack{j = 0\\ \text{$j$ even}}}^{q-2} \rho(\alpha^j+1)\rho(\alpha^{\tau})\rho(\alpha^{-j} + \alpha^{-\tau}) + \sum_{\substack{j = 0\\ \text{$j$ odd}}}^{q-2} \rho(\alpha^j+1)\rho(\alpha^{\tau})\rho(-\alpha^{-j} + \alpha^{-\tau})|\]
\[= |\sum_{j=0}^{q-2} \left(\frac{1 + \rho(\alpha^j)}{2}\right)\rho(\alpha^j+1)\rho(\alpha^{-\tau})\rho(\alpha^{-j})\rho(\alpha^{-\tau}\alpha^j + 1)\] \[+ \sum_{j=0}^{q-2} \left(\frac{1 - \rho(\alpha^j)}{2}\right)\rho(\alpha^j+1)\rho(\alpha^{-\tau})\rho(\alpha^{-j})\rho(\alpha^{-\tau}\alpha^j - 1)|\]
\[\leq \frac{1}{2}(|\sum_{x \in \mathbb{F}_q} \rho(x+1)\rho(x)\rho(x + \alpha^{\tau})|
+ |\sum_{x \in \mathbb{F}_q} \rho(x+1)\rho(x + \alpha^{\tau})|\]
\[|\sum_{x \in \mathbb{F}_q} \rho(x+1)\rho(x)\rho(x - \alpha^{\tau})|
+ |\sum_{x \in \mathbb{F}_q} \rho(x+1)\rho(x - \alpha^{\tau})|) + 1.\]
For $\tau \neq 0$ or $(q-1)/2,$ the Weil bound (Theorem \ref{wb}) implies that the magnitude of the first sum is less than or equal to $2\sqrt{q}+3,$ the magnitude of the second sum is less than or equal to $\sqrt{q} + 2,$ the magnitude of the third sum is less than or equal to $2\sqrt{q}+3,$ and the magnitude of the fourth sum is less than or equal to $\sqrt{q}+2.$ So, if $\tau \neq 0$ or $(q-1)/2,$ then $|\mathcal{C}_{\mathbf{s}, \mathbf{s}[(q-1)/2 - 1]}(\tau)| \leq 3\sqrt{q} + 6.$

If $\tau = 0,$ then the magnitude of the first sum is $0,$ the magnitude of the second sum is $q,$ and by the Weil bound, the magnitude of the third sum is less than or equal to $2 \sqrt{q}+3$ and the magnitude of the fourth sum is less than or equal to $\sqrt{q}+2.$ So, $|\mathcal{C}_{\mathbf{s}, \mathbf{s}[(q-1)/2 - 1]}(0)| \leq \frac{1}{2}(q + 3\sqrt{q}+7).$ Likewise, $|\mathcal{C}_{\mathbf{s}, \mathbf{s}[(q-1)/2 - 1]}((q-1)/2)| \leq \frac{1}{2}(q + 3\sqrt{q}+7).$

By similar arguments, one can show that for $\tau \neq 0$ or $(q-1)/2,$ $|\mathcal{C}_{\mathbf{s}, \mathbf{s}[(q-1)/2 + 1]}(\tau)| \leq 3\sqrt{q} + 6.$ Furthermore, one can also show that $|\mathcal{C}_{\mathbf{s}, \mathbf{s}[(q-1)/2 + 1]}(0)| \leq \frac{1}{2}(q + 3\sqrt{q}+7)$ and $|\mathcal{C}_{\mathbf{s}, \mathbf{s}[(q-1)/2 + 1]}((q-1)/2)| \leq \frac{1}{2}(q + 3\sqrt{q}+7). \qed$

If (one of) $(q-1) \pm 1$ is a multiplier of $S,$ then for some $\tau = 0,...,q-2,$ $\mathcal{C}_{\mathbf{s},\mathbf{s}[(q-1) \pm 1]}(\tau) = q-1.$ However, by Lemma \ref{weakbd}, $\mathcal{C}_{\mathbf{s},\mathbf{s}[(q-1)\pm 1]}(\tau) \leq \text{max}\lbrace 3\sqrt{q}+6, \frac{1}{2}(q + 3\sqrt{q} + 7)\rbrace.$ 
 
Now, $q - 1 \leq 3\sqrt{q} + 6$ if and only if $q-3\sqrt{q} \leq 7.$ Similarly, $q-1 \leq \frac{1}{2}(q + 3\sqrt{q}+7)$ if and only if $q - 3\sqrt{q} \leq 9.$ But, for $x \geq 4,$ $x - 3\sqrt{x}$ is an increasing function, and for $q = 5^2,$ $q - 3\sqrt{q} = 10>7,9.$ Thus, if $q \geq 25,$ then $(q-1)/2 \pm 1$ are not multipliers of $S.$ It can be checked directly that for $q < 25,$ $(q-1)/2\pm 1$ are multipliers of $S$ if and only if $q $ equals $9$. Thus, the following proposition is a direct consequence of Lemma \ref{weakbd}.
 
\begin{proposition} \label{exceptional cases} Let $q$ be a prime power such that $q \equiv 1 \pmod{4},$ and let $S$ be an SLCE almost difference set defined over $\mathbb{F}_q^*.$ Then $(q-1)/2 \pm 1$ are multipliers of $S$ exactly when $q$ equals $9$.
\end{proposition}

These exceptional cases having been dealt with, we now turn to the work of using Theorem \ref{nec} to show that the multiplier group of an SLCE almost difference set over a prime field is trivial.
\begin{notation} Let $S_1 := \lbrace j \in (\mathbb{Z}/(p-1)\mathbb{Z})^*: j^{\prime}<(p-1)/2 \rbrace.$
\end{notation}
\begin{lemma} \label{like Akiyama} Let $p$ be an odd prime, and let $S$ be an SLCE almost difference set over $\mathbb{F}_p^*.$ Let $t \in (\mathbb{Z}/(p-1)\mathbb{Z})^*$ be a multiplier of $S$ not equal to $\pm1.$ Then there exists $a \in S_1$ such that $aS_1 = S_1$ and $a \neq 1.$
\end{lemma}
\emph{Proof.} Since the multipliers of $S$ form a group, $t^{-1}$ is also a multiplier of $S.$ Hence, by Theorem \ref{nec}, $t^{-1}S_0 = S_0$ or $t^{-1}S_0 \cap S_0 = \emptyset.$ 

Assume first that $t^{-1}S_0 = S_0.$ Then for each $j \in (\mathbb{Z}/(p-1)\mathbb{Z})^*,$ $(j^{-1})^{\prime} < (p-1)/2$ if and only if $(tj^{-1})^{\prime} < (p-1)/2.$ Hence, for each $j \in (\mathbb{Z}/(p-1)\mathbb{Z})^*,$ $(j)^{\prime} < (p-1)/2$ if and only if $(tj)^{\prime} < (p-1)/2.$ Therefore, $S_1 = tS_1.$ Furthermore, since $1 \in S_1, $ $S_1 = tS_1$ implies $t \in S_1.$

Now assume that $t^{-1}S_0 \cap S_0 = \emptyset.$ It then follows from the fact that $|-t^{-1}S_0| = |S_0| = |S_0^c|$ and the definition of $S_0$ that $-t^{-1}S_0 = S_0.$ Thus, we can apply the above argument with ``$-t$'' in place of ``$t$'' to deduce that $S_1 = -tS_1$ (and that $-t \in S_1$). $\qed$\\

As it turns out, the problem of deciding which elements $a \in (\mathbb{Z}/(p-1)\mathbb{Z})^*$ ($a \neq 1$) satisfy the equation $aS_1 = S_1$ is very similar to a probelm which arises in a different context.  A Jacobi sum is called \emph{pure} if some positive integral power of it is real. In \cite{A4}, Akiyama determines the conditions under which Jacobi sums of the form $K(\chi)$ defined over $\mathbb{F}_{p^2}$ are pure. Like our work in this paper, Akiyama's work relies on Stickleberger's Theorem: indeed, he shows that a Jacobi sum is pure exactly when a condition holds that is almost identical to the necessary condition for a residue to be a multiplier given in Lemma \ref{like Akiyama}.
\begin{proposition} \label{similarprop} \cite{A4} Let $p$ be an odd prime, let $k|p^2-1,$ and let $\chi$ be a character on $\mathbb{F}_{p^2}$ of order $k.$ Let $R_1 = \lbrace x \in (\mathbb{Z}/k\mathbb{Z})^*: x^{\prime} \in [1,k/2)\cap \mathbb{Z}\rbrace.$ Then $K(\chi)$ is pure if and only if there exists $a \in R_1$ such that $aR_1 = R_1$ and $p \equiv -a \pmod{k}.$
\end{proposition}
As a result of the similarity between our condition and Akiyama's condition, we are able to apply the methods from \cite{A4} almost directly. Akiyama's classification of pure Jacobi sums breaks into a number of cases, as does our proof that the multiplier group is of $S$ is trivial when $d=1.$ We will explicitly give the proof of our result in two special cases in order to show how the ideas from \cite{A4} translate to our context. 

The proof of the following corollary is a straightforward modification of the proof of Lemma $4$ in \cite{A4}.
\begin{corollary} \label{not3mod4} Let $p$ be a prime congruent to $3$ mod $4,$ and let $S$ be an SLCE almost difference set over $\mathbb{F}_p^*.$ Then the multiplier group of $S$ is trivial.
\end{corollary}  
\emph{Proof.} Let $t\neq \pm 1$ be a multiplier of $S.$ Then, by Lemma \ref{like Akiyama}, there exists $a \in S_1$ such that $a \neq 1$ and $aS_1 = S_1.$ 
Pick an integer $i$ such that 
\begin{eqnarray} \label{first} \frac{p-1}{2^{i+1}} < a^{\prime} \end{eqnarray}
and
\begin{eqnarray} \label{second} a^{\prime} \leq \frac{p-1}{2^i}. \end{eqnarray}
Since $a \neq 1,$ it follows from (\ref{second}) that $(p-1)/2^{i+1} \geq 1$ and so that 
\begin{eqnarray} \label{first half} 2^i \leq (p-1)/2. \end{eqnarray}

Since $p \equiv 3 \pmod{4},$ there exists a congruence class $y \in (\mathbb{Z}/(p-1)\mathbb{Z})^*$ containing $(p-1)/2 - 2^i.$ It follows from (\ref{first half}) that $y \in S_1.$ Since $a^{\prime}$ is odd, we have that 
\[a^{\prime}\left(\frac{p-1}{2} - 2^i\right) \equiv \frac{p-1}{2} - 2^ia^{\prime} \pmod{p-1}.\]
By (\ref{first}) and (\ref{second}), 
\[\frac{p-1}{2} < 2^ia^{\prime}\leq p-1.\]
Ergo, $ay \in S_1^c.$ But this contradicts the fact that $aS_1 = S_1.$

So, if $t$ is a multiplier of $S,$ then $t = \pm 1.$ But, by Proposition \ref{minus one}, $-1$ is never a multiplier of $S.$ Hence, the multiplier group of $S$ is trivial. $\qed$\\

The proof of the next result is a straightforward modification of the proof of Lemma $5$ from \cite{A4}.
\begin{corollary} \label{notdiv8} Let $p$ be a prime congruent to $1$ mod $8$ and greater than $14^2+1.$ Let $S$ be an SLCE almost difference set over $\mathbb{F}_p^*.$ Then the multiplier group of $S$ is trivial.
\end{corollary}
\emph{Proof.} Let $t\neq \pm 1, (p-1)/2 \pm 1$ be a multiplier of $S.$ Then, by Lemma \ref{like Akiyama}, there exists $a \in S_1$ such that $a \neq 1$ and $aS_1 = S_1.$ For $c,d \in \mathbb{Z}^{+},$ set $T(c,d) = \lbrace x \in (\mathbb{Z}/(p-1)\mathbb{Z})^* : x^{\prime} \in [c,d)\cap \mathbb{Z} \rbrace,$ and for $j = 1,2,3,4,$ set $T_j = T(((j-1)(p-1)/4,j(p-1)/4).$

Note that there is a residue $y \in S_1$ containing $(p-1)/2 - a^{\prime}.$ Furthermore, since every integer belonging to a congruence class in $S_1$ is odd, the condition $aS_1 = S_1$ is equivalent to the condition $yS_1 = S_1.$ Note that the hypothesis that $t \neq (p-1)/2 \pm 1$ guarantees that $y \neq \pm 1.$ Hence, we may assume $a \in T_1.$ Indeed, let us begin by assuming $a^{\prime} \in [8,(p-1)/4)\cap \mathbb{Z}.$

Let $i$ be the positive integer such that 
\begin{eqnarray} \label{interval} a^{\prime} \in \left[\frac{p-1}{2^{i+2}},\frac{p-1}{2^{i+1}}\right).\end{eqnarray}
Write $p-1 = 2^em,$ where $e\geq 3$ and $m$ is odd. Let $A,$ $B,$ $C,$ and $D$ be elements of $(\mathbb{Z}/(p-1)\mathbb{Z})^*$ containing $(p-1)/2^e + 2^i,$ $(p-1)/2^e + 2^{i+1},$ $(p-1)/2^e + (p-1)/4 + 2^i,$ and $(p-1)/2^e + (p-1)/4 + 2^{i+1},$ respectively. Note that since $a^{\prime} \leq (p-1)/2^{i+1}$ and $a^{\prime} \geq 8,$ we have that $2^{i+4} \leq 2^{i+1}a^{\prime}\leq p-1.$ Hence, $A,B,C,D \in S_1.$

Assume for the sake of contradiction that $aA,$ $aB,$ $aC,$ and $aD$ are all in $S_1.$ First, note that 
\[\left(\frac{p-1}{2^e} + 2^{i+1}\right) - \left(\frac{p-1}{2^e} + 2^i \right) = 2^i,\]
and by (\ref{interval}), 
\[2^ia^{\prime} \in \left[\frac{p-1}{4},\frac{p-1}{2}\right).\]
Hence, $aB-aA \in T_2.$ So, if $aA \in T_2.$ then $aB \in T_3$ and so $aB \in S_1^c,$ which contradicts our assumption. Thus, $aA \in T_1,$ and $aB \in T_2.$

Since $a^{\prime}$ is odd, we consider the following two cases.\\
$1)$ ($a^{\prime} \equiv 1 \pmod{4}$) In this case, 
\[a^{\prime}\left(\frac{p-1}{2^e} + \frac{p-1}{4} + 2^{i+1}\right) \equiv \frac{p-1}{4} + a^{\prime}\left(\frac{p-1}{2^e} + 2^{i+1}\right) \pmod{p-1}.\]
Hence, since $aB \in T_2,$ it follows that $aD \in T_3$ and so $aD \in S_1^c,$ which contradicts our assumption.\\
$2)$ ($a^{\prime} \equiv 3 \pmod{4}$) In this case, 
\[a^{\prime}\left(\frac{p-1}{2^e} + \frac{p-1}{4} + 2^i\right) \equiv \frac{3(p-1)}{4} + a^{\prime}\left(\frac{p-1}{2^e} + 2^i\right) \pmod{p-1}.\]
Hence, since $aA \in T_1,$ it follows that $aC \in T_4$ and so $aC \in S_1^c,$ which contradicts our assumption. Thus, $aA,$ $aB,$ $aC,$ and $aD$ cannot all lie in $S_1$ simultaneously; the equation $aS_1 = S_1$ cannot be true. 

It remains to consider the cases $a^{\prime} = 3,5,7.$ But, note that $aS_1 = S_1$ implies $a^2S_1 = S_1.$ Also, $(p-1)/4 > 7^2$ implies that for each of these choices of $a^{\prime},$ $(a^2)^{\prime} \in [8,(p-1)/4)$ and so the above argument can be applied to obtain a contradiction.

Thus, if $t$ is a multiplier of $S,$ then $t = \pm 1, (p-1)/2 \pm 1.$ But, by Proposition \ref{minus one} and Proposition \ref{similarprop}, $-1$ and $(p-1) \pm 1$ are not multipliers of $S.$ Hence, the multiplier group of $S$ is trivial. $\qed$\\

For proofs of the next three corollaries, see the proofs of Lemma $6,$ Lemma $7,$ and Lemma $8,$ respectively, in \cite{A4}.

\begin{corollary} \label{notprime3} Let $p$ be a prime satisfying $p-1 = 4m$ for some integer $m$ such that $(m,3) = 1.$  Let $p$ be greater than $10^2+1.$ Let $S$ be an SLCE almost difference set over $\mathbb{F}_p^*.$ Then the multiplier group of $S$ is trivial.
\end{corollary}

\begin{corollary} \label{squarefree} Let $p$ be a prime satisfying $p-1 = 4m$ for some integer $m$ which is odd and not square free. Let $p$ be greater than $46^2 + 1.$ Let $S$ be an SLCE almost difference set over $\mathbb{F}_p^*.$ Then the multiplier group of $S$ is trivial.
\end{corollary}

\begin{corollary} \label{notcong3} Let $p$ be a prime satisfying $p-1 = 12m$ for some integer $m$ that has a prime factor greater than $6.$ Further assume that $(m,6) = 1.$ Let $p$ be greater than $70^2 + 1.$ Let $S$ be an SLCE almost difference set over $\mathbb{F}_p^*.$ Then the multiplier group of $S$ is trivial.
\end{corollary}

Let $S$ be an SLCE almost difference set over $\mathbb{F}_p^*.$ Assume, to begin with, that $p>70^2 + 1.$ By Corollary \ref{not3mod4}, the multiplier group of $S$ is trivial unless $4|(p-1).$ So, assume $p-1 = 4m,$ for some integer $m.$ By Corollary \ref{notdiv8}, the multiplier group of $S$ is trivial unless $m$ is odd. Assume $m$ is odd. By Corollary \ref{notprime3}, the multiplier group of $S$ is trivial unless $3|m.$ So, assume $p-1 = 12m^{\prime},$ for some integer $m^{\prime}.$ By Corollaries \ref{notdiv8} and \ref{squarefree}, unless $m^{\prime}$ is square free and relatively prime to $6,$ then the multiplier group of $S$ is trivial. Assume $m^{\prime}$ is indeed square free and relatively prime to $6.$ Then, by Corollary \ref{notcong3}, the multiplier group of $S$ is trivial unless $m^{\prime} = 5,$ in which case $p<70^2 + 1,$ contradicting our assumption.

Thus, when $p>70^2 + 1,$ the multiplier group of $S$ is trivial. Using a simple Python program, we verified that for each $p<70^2 + 1$ such that $p \equiv 1 \pmod{4},$ the multiplier group of $S$ is trivial. For each such prime, we checked every element of $(\mathbb{Z}/(p-1)\mathbb{Z})^*$ not equal to $\pm 1, \pm ((p-1)/2-1)$ against Lemma \ref{like Akiyama} to see whether it could be ruled out as a multiplier. It turns out that Lemma \ref{like Akiyama} rules out almost every potential multiplier. Somewhat oddly, the only exceptions are $t = \pm 11$ and $t= \pm 19$ when $p = 61.$ These exceptional cases are, in fact, quite similar to some exceptional cases noted in Akiyama's paper (see \cite[p. 99]{A4}). One can check directly that $\pm 11$ and $\pm 19$ are not multipliers of $S$ in the case that $p = 61.$ Hence, we have the following theorem.
\begin{theorem} Let $p$ be an odd prime, let $S$ be an SLCE almost difference set over $\mathbb{F}_p^*,$ and let $\mathbf{s}$ be the SLCE sequence corresponding to $S.$ Then the multiplier group of $S$ is trivial, and $\mathcal{F}_1 = \lbrace \mathbf{s}[t^{\prime}] : t \in (\mathbb{Z}/(p-1)\mathbb{Z})^* \rbrace$ is a family of $\phi(p-1)$ shift inequivalent decimations of $\mathbf{s}.$
\end{theorem} 

Our work in this paper suggests several (as of yet open) problems for future research. It would be interesting to extend our results to $M$-ary Sidelnikov sequences (where $M$ is any divisor of $q-1$). It would also be interesting to obtain explicit classifications of the multiplier groups of SLCE almost difference sets when $d>1.$ Finally, it would be interesting to obtain new bounds on the cross-correlation values of a Sidelnikov sequence with one of its decimations.

\section*{Acknowledgements}

The research of \c{S}aban Alaca was supported by a Discovery Grant from the Natural
Sciences and Engineering Research Council of Canada (RGPIN-2015-05208). 
This paper is part of Goldwyn Millar's PhD thesis work; his research was supported by an Ontario Graduate Scholarship.

\end{document}